\documentclass[reqno,12pt]{amsart}
\usepackage{graphicx}
\usepackage{epsfig}
\usepackage{amssymb,amsmath,amsfonts,amstext,mathrsfs}
\usepackage{latexsym}
\usepackage{url}
\usepackage{fancyhdr}
\usepackage{color}
\usepackage[unicode]{hyperref}

\newtheorem{thm}{Theorem}

\newtheorem{prop}{Proposition}
\newtheorem{conj}{Conjecture}
\newtheorem*{prob}{Problem}

\hypersetup{
colorlinks=true,
linkcolor=blue,
citecolor=blue,
urlcolor=blue,
filecolor=blue,
bookmarksnumbered=true,
pdfstartview=FitH,
pdfhighlight=/N
}

\def\d{\,{\rm{d}}}

\def\m{\,{\mathfrak{m}}}

\begin{document}

\title[Fourier coefficients of the $?(x)$ function]{Fourier-Stieltjes coefficients of \\the Minkowski question mark function}

\author[Giedrius Alkauskas]{Giedrius Alkauskas}
\address{Vilnius University, Department of Mathematics and Informatics, Naugarduko 24, LT-03225 Vilnius, Lithuania, \& Institute of Mathematics, Department of Integrative Biology,
Universit\"{a}t f\"{u}r Bodenkultur Wien, Gregor Mendel-Stra{\ss}e 33, A-1180 Wien, Austria,}
\email{giedrius.alkauskas@gmail.com}

\begin{abstract}
In this paper we investigate the Fourier-Stieltjes coefficients of the Minkowski question mark function. In 1943, R. Salem asked whether these coefficients vanish at infinity. We propose the refined conjecture which implies the affirmative answer to Salem's problem. Further, infinite linear identities among these Fourier-Stieltjes coefficients are proved. We also provide a method to numerically calculate special values at integers of the associated zeta function with a high precision (more than 30 digits).
\end{abstract}
\pagestyle{fancy}
\fancyhead{}
\fancyhead[LE]{{\sc Fourier coefficients of $?(x)$}}
\fancyhead[RO]{{\sc G. Alkauskas}}
\fancyhead[CE,CO]{\thepage}
\fancyfoot{}
\date{August 24, 2010. {\it Last update\/}: May 03, 2012.}
\subjclass[2010]{Primary 11A55, 26A30; Secondary 11M41, 42A38}
\keywords{The Minkowski question mark function, Rachjman measure, zeta functions, special values of zeta functions, measure convolution squares}
\thanks{This paper was written in 2010 (the main part) and 2012 (Subsections \ref{1.3} and \ref{2.4}). The author gratefully acknowledges support from the Austrian Science Fund (FWF) under the project Nr. P20847-N18 (in 2010), and support by the Lithuanian Science Council whose postdoctoral fellowship is being funded by European Union Structural Funds project ``Postdoctoral Fellowship Implementation in Lithuania" (in 2012).}
\maketitle

\section{Introduction and main results}
\subsection{Introduction}The Minkowski question mark function $?(x):[0,1]\mapsto[0,1]$, is defined by
\begin{eqnarray}
?([0,a_{1},a_{2},a_{3},\ldots])=2\sum\limits_{i=1}^{\infty}(-1)^{i+1}2^{-\sum_{j=1}^{i}a_{j}},
\label{min}
\end{eqnarray}
where $x=[0,a_{1},a_{2},a_{3},\ldots]$ stands for the representation of $x$ by a regular continued fraction. If $x$ is rational, we can use either finite or infinite expansion; they both give the same value for $?(x)$. This function is continuous, strictly increasing and singular with respect to the Lebesgue measure. It satisfies the functional equations
\begin{eqnarray*}
?(x)=\left\{\begin{array}{c@{\qquad}l} 1-?(1-x),
\\ 2?(\frac{x}{x+1}). \end{array}\right.
\end{eqnarray*}

The \emph{extended} Minkowski question mark function is defined as
$F(x)=?(\frac{x}{x+1})$, $x\in[0,\infty]$. Thus, for $x\in[0,1]$, we
have $?(x)=2F(x)$. In terms of $F(x)$, the above functional
equations can be rewritten as
\begin{eqnarray}
2F(x)=\left\{\begin{array}{c@{\qquad}l} F(x-1)+1 & \mbox{if}\quad x\geq 1,
\\ F(\frac{x}{1-x}) & \mbox{if}\quad 0\leq x<1. \end{array}\right.\label{distr}
\end{eqnarray}
This implies $F(x)+F(1/x)=1$. The overview of available literature on $?(x)$ is contained in \cite{alkauskas_gmj}, and the reader can consult the internet page \cite{mink} for an extensive bibliography list.\\

As was proved by Salem \cite{salem}, the Minkowski question mark function satisfies the Lipschitz condition of order
\begin{eqnarray}
\alpha=\frac{\log 2}{2\log \frac{\sqrt{5}+1}{2}}=0.72021004_{+}.
\label{al}
\end{eqnarray}
Let us, as in \cite{alkauskas_gmj}, define the Laplace-Fourier transform of the Minkowski question mark function by
\begin{eqnarray*}
\m(t)=\int\limits_{0}^{1}e^{xt}\d ?(x).
\end{eqnarray*}
This is the entire function. The symmetry property $?(x)+?(1-x)=1$ implies
\begin{eqnarray}
\m(t)=e^{t}\m(-t).
\label{symm}
\end{eqnarray}
Since $\overline{\m(it)}=\m(-it)$, this shows that $e^{-it/2}\m(it)\in\mathbb{R}$.
Let $d_{n}=\m(2\pi i n)$. Thus, since $e^{-\pi i n}=\pm 1$, we have $d_{n}\in\mathbb{R}$, and so
\begin{eqnarray}
d_{n}=\int\limits_{0}^{1}\cos(2\pi n x)\d ?(x).\label{four}
\end{eqnarray}
\begin{figure}
\begin{center}
\includegraphics[width=260pt,height=380pt,angle=-90]{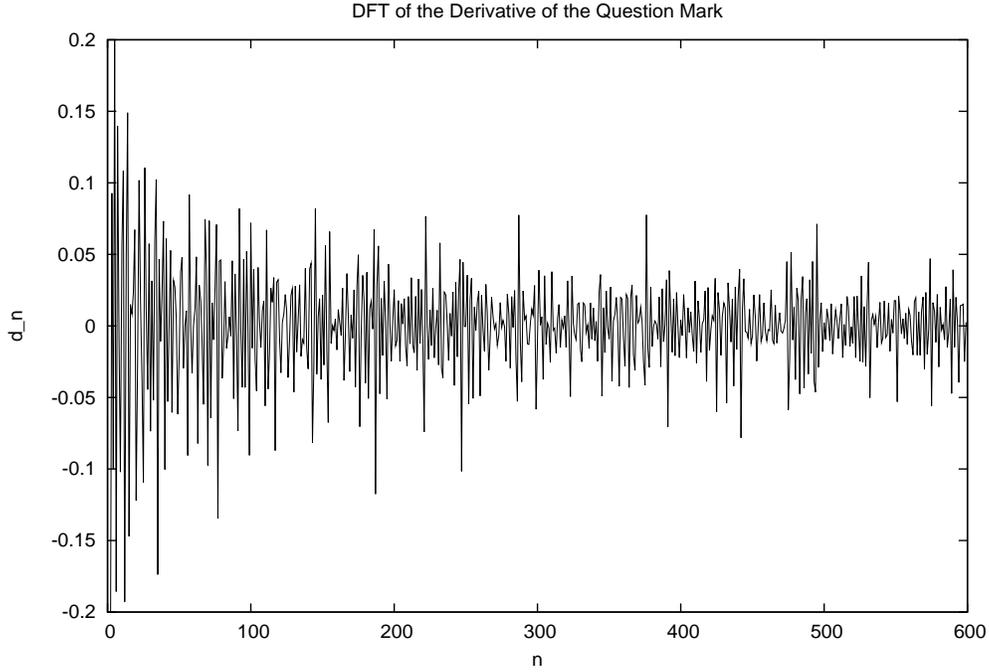}
\caption{The Fourier-Stieltjes coefficients $d_{n}$}
\end{center}
\end{figure}
\begin{figure}
\begin{center}
\includegraphics[width=260pt,height=350pt,angle=-90]{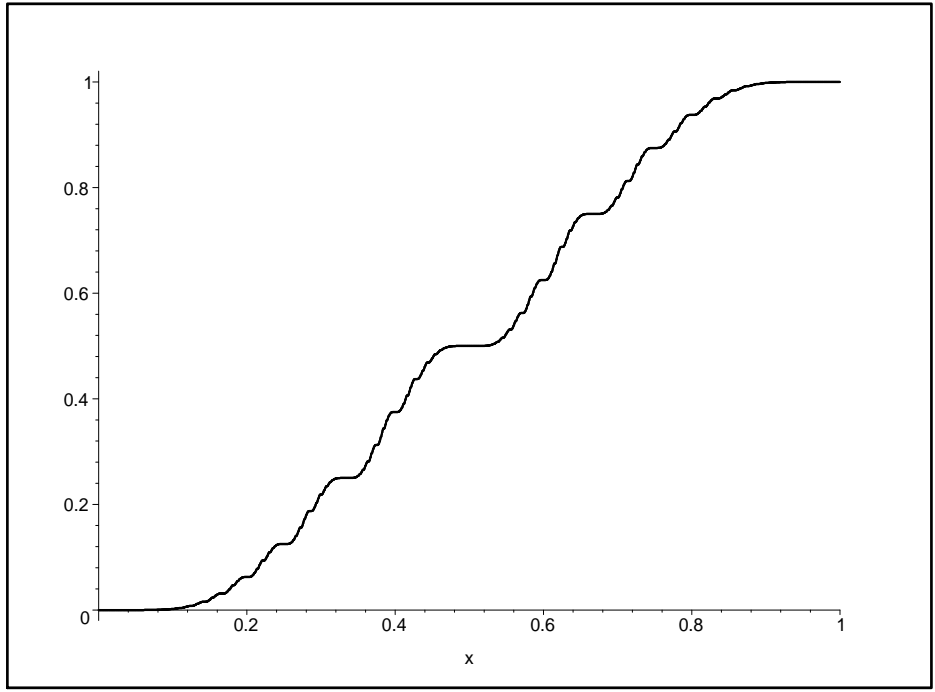}
\caption{The Minkowski question mark function}
\end{center}
\end{figure}
\begin{figure}
\begin{center}
\includegraphics[width=260pt,height=350pt,angle=-90]{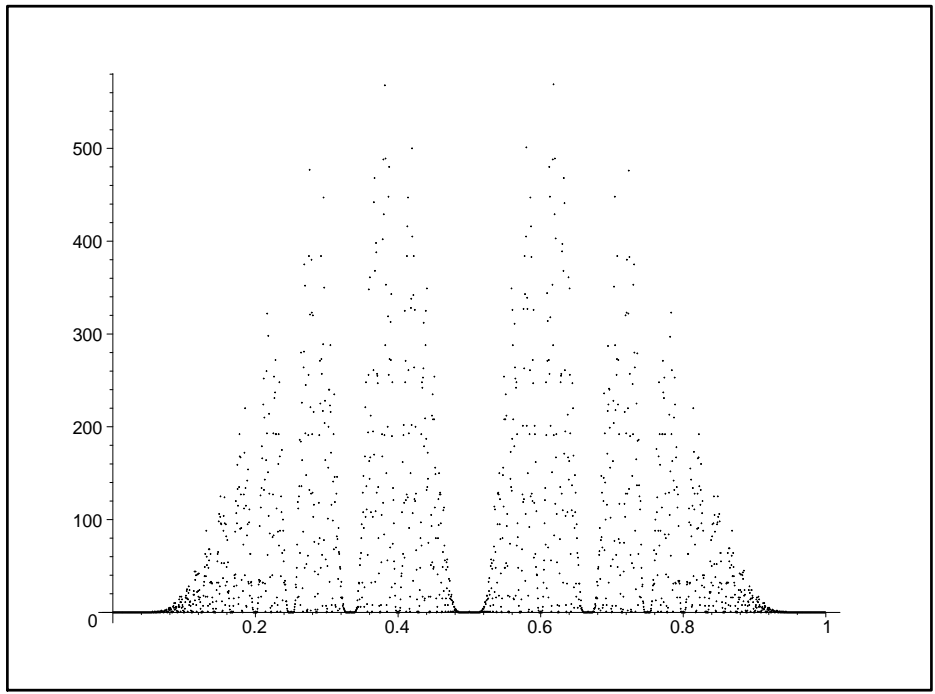}
\caption{The Minkowski measure. More precisely, plot of $N[?(\frac{i+1}{N})-?(\frac{i}{N})]$, $i=0,1,\ldots,N-1$, for large $N$.}
\end{center}
\end{figure}
\begin{figure}
\begin{center}
\includegraphics[width=250pt,height=350pt,angle=-90]{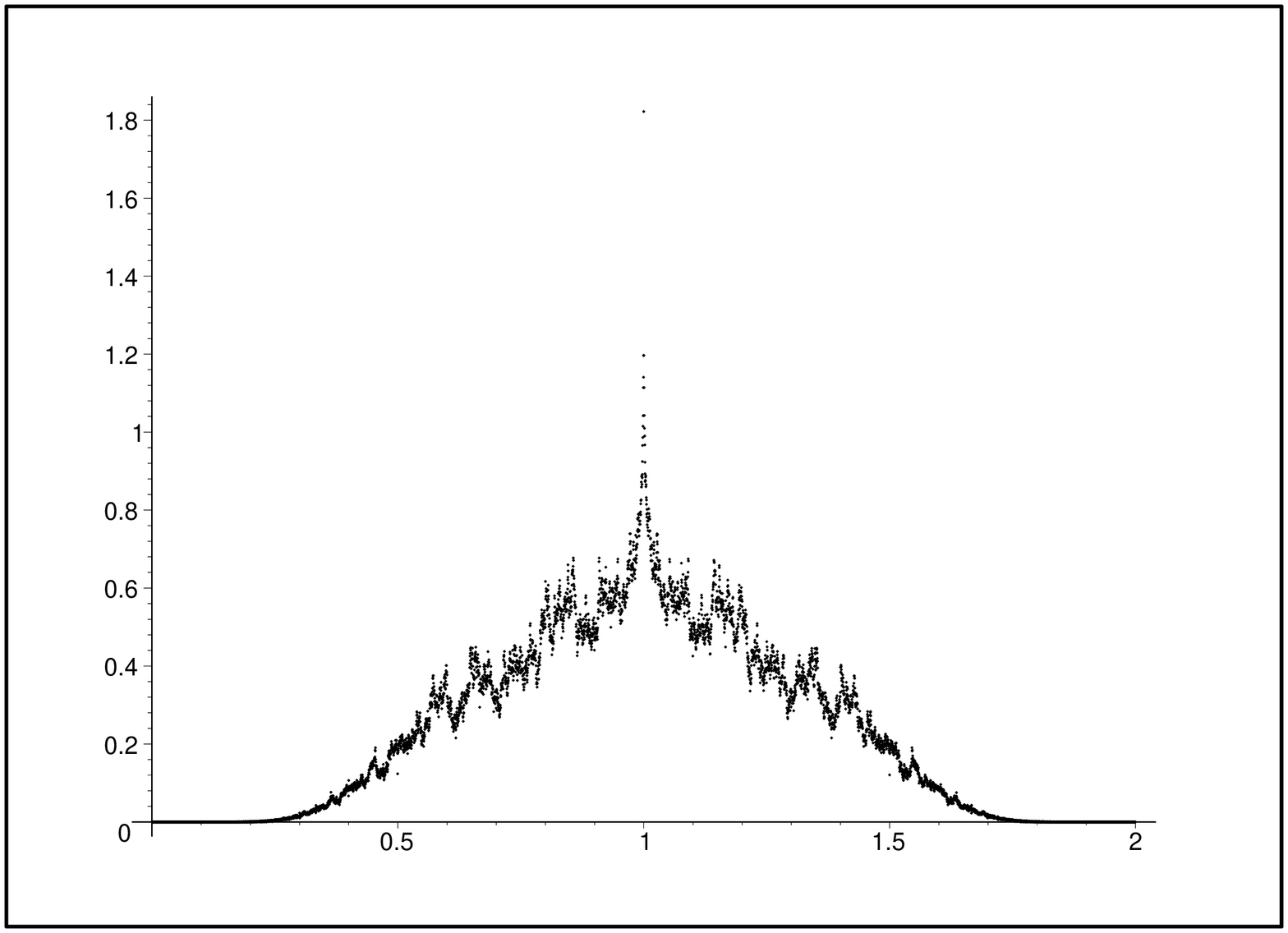}
\caption{Convolution of the Minkowski measure with itself. We conjecture that it is discontinuous only at $x=1$.}
\end{center}
\end{figure}
Note that similar coefficients, defined as $c_{n}=\m(\log2+2\pi i
n)$, are also of importance in the study of the question mark
function and these were investigated in \cite{alkauskas_i}. In \cite{alkauskas_cr}  it was demonstrated that there exist structural linear relations among the constants $d_{n}$, and these are given by
\begin{eqnarray}
d_{m}=\int\limits_{0}^{1}\cos\Big{(}\frac{2\pi m}{x}\Big{)}\d x+
2\sum\limits_{n=1}^{\infty}d_{n}\cdot\int\limits_{0}^{1}\cos(2\pi n x)\cos\Big{(}\frac{2\pi m}{x}\Big{)}\d x,\quad m\in\mathbb{N}.
\label{cr}
\end{eqnarray}
In this paper we prove another collection of more complicated linear relations; the latter arise from the zeta function associated with the coefficents $d_{n}$.
\begin{prop}
For each odd $r\in\mathbb{N}$, there exists a canonical identity of the type
\begin{eqnarray*}
\sum\limits_{n=1}^{\infty}\xi_{r,n}\cdot d_{n}=A_{r}.
\end{eqnarray*}
The series converges absolutely, and the constants $\xi_{r,n}$ and
$A_{r}$ are some explicit constants belonging to the algebra over
$\mathbb{Q}$ generated by Taylor coefficients of $\Gamma(s)$ and
$\zeta(s)$ at $s=1$, also by constants $\pi$, $\frac{1}{\pi}$, $\log
\pi$, $\log n$, $n\in\mathbb{N}$, definite integrals of the form
\begin{eqnarray*}
\int\limits_{n}^{\infty}\Big{[}\sum\limits_{K=1}^{\infty}2^{-K}(\log K)^{i}\sin(2\pi Kx)\Big{]}\cdot\frac{\log^{j}(x)\d x}{x},\quad n\in\mathbb{N},i,j\in\mathbb{N}_{0},
\end{eqnarray*}
and the sum
\begin{eqnarray*}
\sum\limits_{n=1}^{\infty}2^{-n}(\log n)^{r}.
\end{eqnarray*}
\label{thm1}
\end{prop}
The simplest example of this kind of relation arises when $r=1$, and it is given by
\begin{eqnarray}
\frac{3}{4}\cdot\sum\limits_{n=1}^{\infty}2^{-n}\log
n-\frac{1}{4}\cdot\log 2\pi
-\sum\limits_{n=1}^{\infty}\frac{d_{n}}{4n}-\sum\limits_{n=1}^{\infty}\frac{1-d_{n}}{\pi
n}\cdot \int\limits_{n}^{\infty}\frac{\sin(2\pi x)}{5-4\cos(2\pi
x)}\frac{\d x}{x}=0. \label{relation}
\end{eqnarray}
The word ``canonical" in the formulation of the Proposition \ref{thm1} means that the numbering
and the expression for these linear relations does not depend on our choice and, moreover,
they determine the coefficients $d_{n}$ uniquely, provided that the latter arise as a
Fourier-Stieltjes coefficients of a continuous monotone function
$G(x)$ such that $G(0)=0$ and $G(x)+G(1-x)=1$, $x\in[0,1]$.\\

The  function $?(x)$ is a limit distribution of the $n$th generation
of \emph{the Farey tree}. The (permutation of the) latter can be
described by the root $\frac{1}{2}$ and the following rule: each
node $x$ generates two offsprings
\begin{eqnarray*}
x\mapsto \frac{x}{x+1},\quad\frac{1}{x+1}.
\end{eqnarray*}
Thus, for example, let $\mathcal{T}:C[0,1]\mapsto C[0,1]$ be
the operator
\begin{eqnarray*}
\mathcal{T}[f(y)](x)=\frac{1}{2}f\Big{(}\frac{x}{x+1}\Big{)}+\frac{1}{2}f\Big{(}\frac{1}{x+1}\Big{)}.
\end{eqnarray*}
Then in particular we have
\begin{eqnarray*}
d_{n}=\lim\limits_{s\rightarrow\infty}\underbrace{\mathcal{T}\circ\cdots\circ\mathcal{T}}_{s}[\cos(2\pi n y)](x)\text{ for any fixed }x\in[0,1].
\end{eqnarray*}

\subsection{Salem's problem}In 1943, R. Salem \cite{salem} asked whether $d_{n}\rightarrow 0$,
as $n\rightarrow\infty$. Before discussing this deeper, we shortly
overview general properties of Fourier transforms and
Fourier-Stieltjes coefficients of monotone functions. The following
result was proved independently by R. Salem and A. Zygmund.
\begin{thm}[\cite{salemf,zygmund}]
Let $G(x)$ be any non-decreasing function on $[0,1]$, and let
$e_{n}$ and $e(t)$ be the Fourier-Stieltjes coefficients and the
Fourier-Stieltjes transform, defined by
\begin{eqnarray*}
e_{n}=\int\limits_{0}^{1}e^{2\pi i n x}\d G(x),\quad e(t)=\int\limits_{0}^{1}e^{2\pi i xt}\d G(x).
\end{eqnarray*}
Then, if $e_{n}=o(1)$ for $n\rightarrow\infty$, we have also $e(t)=o(1)$ for $t\rightarrow\infty$. (The converse
is, of course, trivial.)
\end{thm}
The question to determine whether a given measure is a
\emph{Rajchman measure} (that is, whose Fourier transform vanishes
at infinity), as far as measures arising from singular monotone
functions are concerned, is a very delicate question. There are
various examples of all sorts. Let, for example, $G(x)$ by the
classical \emph{Cantor ``middle-third" distribution}. Then Hardy and
Ragosinski observed that using a self-similarity property of $G(x)$,
one gets
\begin{eqnarray*}
e_{3n}=\int\limits_{0}^{1}e^{6\pi i nx}\d G(x)=\int\limits_{0}^{1/3}e^{6\pi i nx}\d G(x)+\int\limits_{2/3}^{1}e^{6\pi i nx}\d G(x)=
\Big{(}\frac{1}{2}+\frac{1}{2}\Big{)}\int\limits_{0}^{1}e^{2\pi i nx}\d G(x).
\end{eqnarray*}
So, for $N\in\mathbb{N}$, $e_{3^{N}}=e_{1}$ (which can be easily seen to bee non-zero), and thus
$e_{n}\neq o(1)$. Of course, $G(x)$ is not strictly increasing. But this is not an obstacle. On can construct \cite{salem,reese}
further examples of singular functions which are strictly increasing and whose Fourier-Stieltjes coefficients
still do not vanish at infinity. On the other hand, there exist
singular distributions whose coefficients $e_{n}$ vanish,
and the first example was given by Menchoff in 1916. Due to contributions by Wiener and Wintner \cite{wiener},
Salem \cite{salem2,salem3}, Schaeffer \cite{schaeffer}, Iva\v{s}ev-Musatov \cite{musatov},
we know that for every $\epsilon>0$ there exists a singular monotone distribution $G(x)$ such
that $e_{n}=O(n^{-1/2+\epsilon})$. Note however that if $|e(t)|=O(|t|^{-1-\epsilon})$
for a certain $\epsilon>0$, then $G(x)$ is necessarily absolutely continuous. See also \cite{bluhm,hartman}.\\
\indent Thus, the theorem of Salem and Zygmund tells \emph{a priori}
that the fact $d_{n}=o(1)$ implies $\m(it)=o(1)$. As was noted in \cite{salem}, the general
theorem of Wiener \cite{zygmund} about the Fourier-Stieltjes
coefficients of continuous monotone functions with known modulus of
continuity implies
\begin{eqnarray*}
\sum\limits_{i=1}^{n}\vert d_{i}\vert^{2}=O(n^{1-\alpha})=O(n^{0.2797_{+}}).
\end{eqnarray*}
Here $\alpha$ is given by (\ref{al}). This shows that $|d_{n}|\ll
n^{-0.3601}$ on average. For future purpose we note that, as also
observed in \cite{salem}, the Cauchy-Schwartz inequality implies
\begin{eqnarray}
\Omega(n):=\sum\limits_{i=1}^{n}|d_{i}|=O(n^{1-\alpha/2}). \label{daline}
\end{eqnarray}

Salem's question can be restated in terms of the Farey tree, which we have already defined earlier. Let $\mathcal{T}_{m}$ denote the $m$th generation of this
tree, $m\in\mathbb{N}_{0}$. Thus, $\mathcal{T}_{0}$ consists of a
single rational number $\frac{1}{2}$.
\begin{prob}[R. Salem] Prove or disprove that
\begin{eqnarray*}
\lim\limits_{n\rightarrow\infty}\lim\limits_{m\rightarrow\infty}\frac{1}{2^{m}}\sum\limits_{\frac{p}{q}\in\mathcal{T}_{m}}
e^{2\pi in\frac{p}{q}}=0.
\end{eqnarray*}
\end{prob}
\subsection{Refined conjecture}
\label{1.3}
In \cite{alkauskas_cr} two results related to the above were proved. The next approach towards the solution of Salem's problem we propose jointly with H\aa kan Hedenmalm and Alfonso Montes-Rodriguez.\\
\indent As can be proved, the fact that a Fourier transform of a certain
measure $(\d\mu,[0,1])$ vanishes at infinity follows from the fact
that there exists $k\in\mathbb{N}$ such that
$\underbrace{\mu'\circ\cdots\circ\mu'}_{k}\in L^{1}[0,1]$. This follows from the Lebesgue-Riemann lemma. See \cite{bisbas, gupta, hewitt, saeki} for
related results and problems. In our case it is most likely that $k=2$ is
already sufficient. The  Figures 2-4 show the Minkowski question mark function, the (formal) Minkowski measure $?'(x)$, and the convolution $?'(x)\circ ?'(x)$. The latter picture suggests that $?'(x)\circ ?'(x)$ is a continuous function on $[0,1)\cup(1,2]$, and it has an integrable singularity at $x=1$. To prove these facts, one needs to perform analysis on the the sums of the form
\begin{eqnarray*}
\Sigma(N)=\sum\limits_{i=0}^{N-1}\Big{[}?\Big{(}\frac{i+1}{N}\Big{)}-?\Big{(}\frac{i}{N}\Big{)}\Big{]}^2,\quad
N\in\mathbb{N}.
\end{eqnarray*}
Though, of course, we have an explicit expression for the value of
the question mark function at rational argument, i.e. the formula (\ref{min}), the problem to find asymptotics, or at least a good estimate for the above sum is quite involved. Of course, the Lipschitz property gives
\begin{eqnarray*}
\Sigma(N)\ll N^{1-2\alpha}=N^{-0.4404_{+}},
\end{eqnarray*}
but this is far too rough; we would expect something of order $N^{-1}\log N$. In what follows, we assume $?(x)=0$ for $x\leq 0$ and $?(x)=1$ for $x\geq 1$. We raise the following
\begin{conj}
For every $\beta\in[0,2]\setminus\{1\}$, there exists the constant $C(\beta)\geq 0$ such that
\begin{eqnarray*}
\sum\limits_{i=0}^{N-1}\Big{[}?\Big{(}\frac{i+1}{N}\Big{)}-?\Big{(}\frac{i}{N}\Big{)}\Big{]}\cdot
\Big{[}?\Big{(}\beta-\frac{i}{N}\Big{)}-?\Big{(}\beta-\frac{i+1}{N}\Big{)}\Big{]}\sim\frac{C(\beta)}{N},\text{ as } N\rightarrow\infty,\quad N\in\mathbb{N}.
\end{eqnarray*}
 The function $C(\beta)$ is discontinuous only at $\beta=1$, where  both side limits are $\infty$, but it is an integrable singularity.
\end{conj}
The graph of the function $C(\beta)=C(2-\beta)$ resembles the graph of the convolution $?'(x)\circ ?'(x)$, given in the Figure 4. This conjecture, if proved true, implies the positive answer to Salem's problem and the bound $\m(it)=O(t^{-1/4})$ as $t\rightarrow\infty$. Our common project with H. Hedenmalm, A. Montes-Rodriguez, and N. Moshchevitin is in progress.

\section{The Fourier-Stieltjes coefficients and the zeta function}
\label{section2}
\subsection{Special values at even positive integers}
\label{2.1}
 In this section we provide the method to numerically calculate the constants $d_{n}$ to high accuracy and also describe (infinite) linear relations among these coefficients. Let us define
\begin{eqnarray*}
\mathfrak{M}(s)=\sum\limits_{n=1}^{\infty}\frac{d_{n}}{n^{s}}.
\end{eqnarray*}
Since
\begin{eqnarray*}
\Big{|}\sum\limits_{n=1}^{N}d_{n}\Big{|}=\Bigg{|}\int\limits_{0}^{1}\frac{1-e^{2\pi
i (N+1)x}}{1-e^{2\pi i x}}\d ?(x)\Bigg{|} \leq
\int\limits_{0}^{1}\frac{2\d ?(x)}{|1-e^{2\pi i x}|}<+\infty,
\end{eqnarray*}
the Cauchy-Abel principle and properties of general Dirichlet series imply that the series for $\mathfrak{M}(s)$ converges relatively for $\Re(s)>0$. Note that if $d_{n}$ vanish at infinity, then Salem's lemma \cite{salemf} shows that this series converges for $s=0$ to the value $-\frac{1}{2}$. Further, Abel's identity implies that for $\Re(s)=\sigma$,
\begin{eqnarray*}
\sum\limits_{n=1}^{N}\frac{|d_{n}|}{n^{\sigma}}=
\sum\limits_{n=1}^{N-1}\Omega(n)\Big{(}\frac{1}{n^{\sigma}}-\frac{1}{(n+1)^{\sigma}}\Big{)}+\frac{\Omega(N)}{N^{\sigma}},
\end{eqnarray*}
where $\Omega(n)$ is defined by (\ref{daline}). Thus, from (\ref{daline}) we inherit that the Dirichlet series $\mathfrak{M}(s)$ converges absolutely at least for $\Re(s)>1-\frac{\alpha}{2}=0.6398_{+}$. Let
\begin{eqnarray*}
m_{L}=\int\limits_{0}^{1}x^{L}\d ?(x),\quad L\in\mathbb{N}_{0}.
\end{eqnarray*}
(See \cite{alkauskas_mc,alkauskas_rj} for a certain advance to understand the intrinsic structure of these moments.) We have the standard Fourier series for $x^{L}$ in the interval $[0,1]$:
\begin{eqnarray*}
x^{L}\sim\frac{1}{L+1}-\sum\limits_{n\in\mathbb{Z}\atop n\neq 0}e^{2\pi i nx}\Big{(}\frac{1}{2\pi in}+\frac{L}{(2\pi i n)^2}+\cdots+\frac{L(L-1)\cdots 2}{(2\pi i n)^{L}}\Big{)}.
\end{eqnarray*}
Thus, taking the real part, minding (\ref{four}) and the fact that $\int_{0}^{1}\sin(2\pi n x)\d ?(x)$ vanishes, we get
\begin{eqnarray}
m_{L}=\frac{1}{L+1}-2\sum\limits_{n=1}^{\infty}d_{n}\sum\limits_{2\leq 2v\leq L}(-1)^{v}\frac{L!}{(L-2v+1)!(2\pi n)^{2v}}.
\label{mom}
\end{eqnarray}
And so,
\begin{eqnarray*}
m_{L}=\frac{1}{L+1}-2\sum\limits_{2\leq 2v\leq L}\frac{L!(-1)^{v}}{(L-2v+1)!}\cdot\frac{\mathfrak{M}(2v)}{(2\pi)^{2v}},\quad L\geq 0.
\end{eqnarray*}
(For $L=0$ the empty sum is $0$ by convention). Multiply this by $t^{L}/L!$ and sum over $L\geq 0$. We get
\begin{eqnarray*}
\m(t)=\frac{e^{t}-1}{t}\cdot\Big{(}1-2\sum\limits_{2v=2}^{\infty}\frac{\mathfrak{M}(2v)(-1)^{v}}{(2\pi)^{2v}}\,t^{2v}\Big{)}.
\end{eqnarray*}
Therefore,
\begin{eqnarray*}
\frac{1}{2}-\frac{1}{2}\m(t)\sum\limits_{k=0}^{\infty}B_{k}\frac{t^{k}}{k!}=\sum\limits_{2v=2}^{\infty}\frac{\mathfrak{M}(2v)(-1)^{v}}{(2\pi)^{2v}}\,t^{2v}.
\end{eqnarray*}
Here $B_{k}$ are the classical Bernoulli numbers, given by
\begin{eqnarray*}
\frac{t}{e^{t}-1}=\sum\limits_{k=0}^{\infty}B_{k}\frac{t^{k}}{k!}.
\end{eqnarray*}
This finally gives us the expression of special values of the zeta function $\mathfrak{M}(s)$ in terms of the moments:
\begin{eqnarray}
\frac{\mathfrak{M}(2v)}{(2\pi)^{2v}}=\frac{(-1)^{v+1}}{2}\sum\limits_{L=0}^{2v}\frac{B_{2v-L}m_{L}}{(2v-L)!L!},\quad 2v\geq 2.
\label{evenmo}
\end{eqnarray}
We will see soon that this is also valid for $2v=0$. Of course, the last formula is not using any specific information about the Minkowski question mark function: all reasoning remains valid for the moments and a zeta function associated with Fourier-Stieltjes coefficients of any monotone continuous function $G(x)$ such that $G(0)=0$, $G(x)+G(1-x)=1$. Nevertheless, we know a high-precision method to calculate the moments $m_{L}$ \cite{alkauskas_mc}. This gives us a high precision results for the special values $\mathfrak{M}(2v)$, $2v\geq 2$. Table 1 lists the first few.
\begin{table}
\begin{tabular}{|r|r||r|r|}
\hline
$v$ & $\mathfrak{M}(2v)$ & $v$ & $\mathfrak{M}(2v)$\\
\hline
$1$ & $-0.4185389015363278$ & $6$ & $-0.3699304357431478$\\
$2$ & $-0.3833363407612589$ & $7$ & $-0.3698882692560897$\\
$3$ & $-0.3733723986086854$ & $8$ & $-0.3698777069802058$\\
$4$ & $-0.3707638984421253$ & $9$ & $-0.3698750640759581$\\
$5$ & $-0.3700983784501058$ & $10$ & $-0.3698744030876739$\\
\hline
\end{tabular}
\caption{Sequence $\mathfrak{M}(2v)$}
\end{table}
We also discover that apart from $m_{2L}$, $L\geq 1$, the constants
$\mathfrak{M}(2v+1)$, $v\geq 0$, also contain certain meaningful
information about the Minkowski question mark function. From the
Table 1 and its extension we also get information about the
coefficients $d_{n}$. For example,
\begin{eqnarray*}
d_{1}=\lim\limits_{v\rightarrow\infty}\mathfrak{M}(2v)=-0.36987418271425_{+}.
\end{eqnarray*}
The given digits are in fact all exact digits which this method and
a standard home computer can provide. Nevertheless, this value is
far better then the Riemann-Stieltjes sum based on (\ref{four}) can
give. Unfortunately, The formula
\begin{eqnarray*}
d_{2}=\lim\limits_{v\rightarrow\infty}2^{2v}\cdot\Big{(}\mathfrak{M}(2v)-d_{1}\Big{)}
\end{eqnarray*}
gives only a few correct digits of $d_{2}$. To calculate other
coefficients $d_{n}$ the following method is much better. Simply, if
we expand $\cos(2\pi nx)$ as Taylor series and use (\ref{four}), we
get
\begin{eqnarray*}
d_{n}=\sum\limits_{k=0}^{\infty}\frac{(2\pi n)^{2k}(-1)^{k}}{(2k)!}\,m_{2k}.
\end{eqnarray*}
This is, obviously, the inverse to the system (\ref{mom}). Based on
the method described in \cite{alkauskas_mc}, we get very good values
of $m_{L}$ for $L<130$. Thus, the last formula can give us good
values for the coefficients $d_{n}$ for $2\pi n<65$, or $n\leq 10$.
The results are provided in the Table 2. To finish the numerical part, we remark that to calculate $d_{n}$ for any $n$ with a precision of $4-5$ decimal digits, the best way (yet) is simply to use Riemann-Stieltjes
sums for the integral (\ref{four}).
\begin{table}
\begin{tabular}{|r|r||r|r|}
\hline
$n$ & $d_{n}$ & $n$ & $d_{n}$\\
\hline
$1$ & $-0.36987418271425589511$ & $6$ & $-0.18571787696613298999$\\
$2$ & $-0.23110608380419115403$ & $7$ & $+0.13977897406302915392$\\
$3$ & $+0.09276672356657101657$ & $8$ & $-0.00611936309545097758$\\
$4$ & $-0.09983104428383632687$ & $9$ & $-0.10205760334150128491$\\
$5$ & $+0.20114256044594273585$ & $10$ &$+0.05670950402333429033$\\
\hline
\end{tabular}
\caption{Sequence $d_{n}$}
\end{table}
\subsection{Preliminary calculations}
We will now show that there exist a countable number of linear relations among $d_{n}$ arising form the zeta function $\mathfrak{M}(s)$. Here ``linear'' means that each constant $d_{n}$ has a factor which is expressed in terms of well-known classical or explicit constants. Here is a very simple linear relation:
\begin{eqnarray*}
\sum\limits_{n=1}^{\infty}d_{n}=-\frac{1}{2}.
\end{eqnarray*}
It holds, provided the answer to Salem's problem is affirmative. Unfortunately, we cannot consider this as a true relation because the series is only conditionally convergent. The ``true'' linear relations are considerably more complicated. As mentioned, the simplest of these is given by (\ref{relation}), being the first case described in Proposition \ref{thm1}, which we will prove in this and the next subsection.\\

Let $\Delta(x)=1-F(x)$, $x\in\mathbb{R}_{+}$. Then (\ref{distr}) shows that $\Delta(x+1)=\frac{1}{2}\Delta(x)$. Extend $\Delta(x)$ to the half-line $(-\infty,0)$ using the last identity. Suppose, $n\in\mathbb{N}$ and $x\in(n-1,n)$. Thus,
\begin{eqnarray*}
\Delta(-x)&=&2^n\Delta(n-x)=2^n(1-F(n-x))=2^n\Big{(}\frac{1}{2}+F(1-n+x)\Big{)}\\
&=&3\cdot2^{n-1}-2^n\Delta(1-n+x)=
3\cdot2^{n-1}-2^{2n-1}\Delta(x).
\end{eqnarray*}
This gives
\begin{eqnarray*}
2^{-n}\Delta(-x)+2^{n-1}\Delta(x)=\frac{3}{2}.
\end{eqnarray*}
Let us introduce
\begin{eqnarray*}
\Phi(x)=2^{\lfloor x\rfloor}\Delta(x).
\end{eqnarray*}
(Note the difference between this function and $\Psi(x)$ which was introduced in \cite{alkauskas_i}). Thus,
\begin{eqnarray*}
\Phi(x)+\Phi(-x)=\frac{3}{2},\text{ if }x\in
\mathbb{R}\setminus\mathbb{Z},\quad \Phi(x)+\Phi(-x)=2\text{ for
}x\in\mathbb{Z}.
\end{eqnarray*}
Additionally, for $x\in\mathbb{R}$, $\Phi(x+1)=\Phi(x)$.
Consequently, the function $\Phi(x)$ is continuous everywhere except
at integer points, and it has the corresponding Fourier expansion
\begin{eqnarray}
\Phi(x)\sim\frac{3}{4}+\sum\limits_{n=1}^{\infty}\frac{\widehat{d}_{n}}{2\pi n}\cdot\sin(2\pi nx).
\label{ser}
\end{eqnarray}
The coefficients $d_{n}$ and $\widehat{d}_{n}$, $n\in\mathbb{N}$, are related by
\begin{eqnarray}
d_{n}=2\int\limits_{0}^{1}\cos(2\pi nx)\d (F(x)-1)=1-4\pi
n\int\limits_{0}^{1}(1-F(x))\sin(2\pi nx)\d x=1-\widehat{d}_{n}.
\label{perein}
\end{eqnarray}
Note that we get a linear relation among the coefficients $d_{n}$
specializing the series (\ref{ser}). Since $2F(x)=?(x)$ and since
for $x\in(0,1)$, $\sum_{n=1}^{\infty}\sin(2\pi n x)(2\pi
n)^{-1}=1/4-x/2$, we get the identity
\begin{eqnarray}
?(x)-x=\sum\limits_{n=1}^{\infty}\frac{d_{n}}{\pi n}\cdot\sin(2\pi nx),\quad x\in[0,1].
\label{trig}
\end{eqnarray}
 Indeed, $?(x)-x$ is a continuous function, and $?(0)-0=?(1)-1$. Further, as we have seen in the beginning of Subsection \ref{2.1},
 the series on the right is absolutely and uniformly convergent to a continuous function, and it converges to $?(x)-x$ in $L^{2}[0,1]$. Hence, we have an identity pointwise. The graph of the function
 $?(x)-x$ can be found in \cite{wiki}. For $x=\frac{1}{3}$, $x=\frac{1}{4}$ and $x=\frac{1}{6}$, the identity (\ref{trig}) gives
\begin{eqnarray*}
\sum\limits_{n=0}^{\infty}\Big{(}\frac{d_{3n+1}}{3n+1}-\frac{d_{3n+2}}{3n+2}\Big{)}=-\frac{\pi}{6\sqrt{3}};& &
\sum\limits_{n=0}^{\infty}(-1)^{n}\frac{d_{2n+1}}{2n+1}=-\frac{\pi}{8};\\
\quad\sum\limits_{n=0}^{\infty}(-1)^{n}\Big{(}\frac{d_{3n+1}}{3n+1}+\frac{d_{3n+2}}{3n+2}\Big{)}&=&-\frac{13\pi}{48\sqrt{3}}.
\end{eqnarray*}
Of course, these three identities are not specific to the Minkowski question mark function: they are using only the symmetry property $?(x)+?(1-x)=1$ and the values $?(\frac{1}{n})=2^{1-n}$, $n=3,4,6$. Thus, to identify $?(x)$ uniquely among symmetric distributions we must choose a dense countable set $S\subset(0,\frac{1}{2})$, and present an identity for each real number in a set $S$. This is not practical and not canonical, and we rather proceed in a different way, which, most importantly, provides a canonical collection of linear relations, very different from a much simpler and also canonical system, proved in \cite{alkauskas_cr}; that is, (\ref{cr}). 
\subsection{Properties of the zeta function}
Let
\begin{eqnarray*}
\mathfrak{L}(s):=\sum\limits_{n=1}^{\infty}\frac{\widehat{d}_{n}}{n^{s}}.
\end{eqnarray*}
Then (\ref{perein}) implies $\mathfrak{M}(s)=\zeta(s)-\mathfrak{L}(s)$; here $\zeta(s)$ is the Riemann zeta function. Consequently, $\mathfrak{L}(s)$ is meromorphic for $\Re(s)>0$, has a simple pole at $s=1$ with a residue $1$. Let $0<\Re(s)<1$. Now, consider the following integral \cite{alkauskas_i}. 
\begin{eqnarray*}
& &\zeta_{\mathcal{M}}(s)\Gamma(s+1)\\
&:=&\int\limits_{0}^{\infty}x^{s}\d
(F(x)-1)=\int\limits_{0}^{\infty}sx^{s-1}(1-F(x))\d x\\
&=&\int\limits_{0}^{\infty}sx^{s-1}2^{-\lfloor x\rfloor}\sum\limits_{n=1}^{\infty}\frac{\widehat{d}_{n}}{2\pi n}\cdot\sin(2\pi
nx)\d x+\frac{3}{4}\int\limits_{0}^{\infty}sx^{s-1}2^{-\lfloor x\rfloor}\d x\\
&=&\frac{3}{4}\sum\limits_{n=1}^{\infty}2^{-n}n^{s}+\sum\limits_{K=1}^{\infty}\sum\limits_{n=1}^{\infty}2^{-K}\frac{\widehat{d}_{n}}{2\pi  n}
\int\limits_{0}^{K}sx^{s-1}\sin(2\pi nx)\d x\\
\end{eqnarray*}
\begin{eqnarray}
&=&\frac{3}{4}\sum\limits_{n=1}^{\infty}2^{-n}n^{s}+\sum\limits_{K=1}^{\infty}
\sum\limits_{n=1}^{\infty}2^{-K}\frac{\widehat{d}_{n}}{2\pi n^{s+1}}\int\limits_{0}^{nK}sy^{s-1}\sin(2\pi
y)\d y \label{funk1}\\
&=&\frac{3}{4}\sum\limits_{n=1}^{\infty}2^{-n}n^{s}-(2\pi)^{-(s+1)}\Gamma(s+1)\cos\frac{\pi
(s+1)}{2}\cdot\mathfrak{L}(s+1)\nonumber\\
&-&\sum\limits_{K=1}^{\infty}\sum\limits_{n=1}^{\infty}2^{-K}\frac{\widehat{d}_{n}}{2\pi n^{s+1}}
\int\limits_{nK}^{\infty}sy^{s-1}\sin(2\pi
y)\d y.\label{funk2}
\end{eqnarray}
Here we used the identities
\begin{eqnarray*}
\int\limits_{0}^{\infty}2^{-\lfloor x\rfloor}g(x)\d x&=&\sum\limits_{K=1}^{\infty}2^{-K}\int\limits_{0}^{K}g(x)\d x;\\
\int\limits_{0}^{\infty}y^{s-1}\sin(2\pi y)\d y&=&-\Gamma(s)(2\pi)^{-s}\cos\frac{\pi(s+1)}{2},\quad -1<\Re(s)<1.
\end{eqnarray*}
We know that $\zeta_{\mathcal{M}}(s)\Gamma(s+1)$  is the entire function. Moreover, the formula (\ref{funk1}) is valid for $\Re(s)>0$, and the formula (\ref{funk2}) is valid for $0<\Re(s)<1$. However, the last sum in (\ref{funk2}) is an analytic function for $\Re(s)<1$. Indeed,
\begin{eqnarray*}
\Big{|}\int\limits_{nK}^{\infty}sy^{s-1}\sin(2\pi y)\d y\Big{|}\ll |s|(nK)^{\Re(s)-1}
\end{eqnarray*}
after integrating once by parts.
So, the formula (\ref{funk2}) gives an analytic continuation of
\begin{eqnarray*}
\Xi(s)=(2\pi)^{-s}\Gamma(s)\cos\frac{\pi s}{2}\cdot\mathfrak{L}(s)
\end{eqnarray*}
to the half plane $\Re(s)<2$. Thus, $\Xi(s)$ is an entire function. Standard analysis then shows that $\mathfrak{L}(s)$ is meromorphic function with a single simple pole at $s=1$, and $\mathfrak{L}(-k)=0$ for $k=0,2,4,\ldots$ Since $\zeta(-k)=0$ for $k=2,4,6,\ldots$ and $\zeta(0)=-\frac{1}{2}$, we have proved the following
\begin{prop}The zeta function $\mathfrak{M}(s)$ extends as an entire function to the whole complex plane with trivial zeros $\mathfrak{M}(-k)=0$ for
$k=2,4,6,\ldots$, and $\mathfrak{M}(0)=-\frac{1}{2}$.
\end{prop}
Let $\overline{f}$ be a monotone continuous function,
$\overline{f}(0)=0$, $\overline{f}(x)+\overline{f}(1-x)=1$. Lets us
define $\overline{F}(x+n)=1-2^{-n}+2^{-n}\overline{f}(x)$ for
$x\in[0,1]$, $n\in\mathbb{N}_{0}$. The equalities (\ref{funk1}) and
(\ref{funk2}) also do not use any specific information about the Minkowski
question mark function: they also hold if $1-F(x)$ is replaced by a
generic $1-\overline{F}(x)$ (of course, to show that
$\zeta_{\mathcal{M}}(s)\Gamma(s+1)$ is entire we are explicitly
using a fact that $\d F(x)$ ``kills" every power $x^{-L}$ at $x=0$).
Nevertheless, the Minkowski question mark function is characterised
by the invariance of the l.h.s. of (\ref{funk1}) and (\ref{funk2})
under the map $s\mapsto -s$ \cite{alkauskas_i}. Thus, the function
$\zeta_{\mathcal{M}}(s)\Gamma(s+1)$ is even function, hence its odd
order derivatives at $s=0$ vanish. In particular, the first
derivative of (\ref{funk2}) at $s=0$ vanishes. Thus,
\begin{eqnarray}
0&=&\frac{3}{4}\sum\limits_{n=1}^{\infty}2^{-n}\log n-\frac{d}{\d s}\Big{(}(2\pi)^{-(s+1)}\Gamma(s+1)\cos\frac{\pi
(s+1)}{2}\cdot\mathfrak{L}(s+1)\Big{)}\Big{|}_{s=0}\label{ko}\\
&-&\sum\limits_{K=1}^{\infty}
\sum\limits_{n=1}^{\infty}2^{-K}\frac{\widehat{d}_{n}}{2\pi n}\int\limits_{nK}^{\infty}y^{-1}\sin(2\pi
y)\d y.
\nonumber
\end{eqnarray}
Note that $\Gamma'(1)=-\gamma$, where $\gamma$ is the Euler-Mascheroni constant. Further,
\begin{eqnarray*}
\zeta(s+1)\cos\frac{\pi(s+1)}{2}=\Big{(}\frac{1}{s}+\gamma+\cdots\Big{)}\Big{(}-\frac{\pi}{2}s+\frac{\pi^{3}}{48} s^{3}+\cdots\Big{)}=-\frac{\pi}{2}-\frac{\pi\gamma}{2}s+\cdots.
\end{eqnarray*}
Thus, the middle term on the r.h.s. of (\ref{ko}) is equal to
\begin{eqnarray*}
-\frac{1}{4}\log 2\pi-\frac{1}{4}\cdot\mathfrak{M}(1).
\end{eqnarray*}
Further, recall that $\widehat{d}_{n}=1-d_{n}$. Note that
\begin{eqnarray*}
\sum\limits_{K=1}^{\infty}2^{-K}\int\limits_{nK}^{\infty}y^{-1}\sin(2\pi
y)\d y=\sum\limits_{K=1}^{\infty}2^{-K}\int\limits_{n}^{\infty}y^{-1}\sin(2\pi
Ky)\d y.
\end{eqnarray*}
Finally,
\begin{eqnarray*}
\sum\limits_{K=1}^{\infty}2^{-K}\sin(2\pi Ky)=\Im(\sum\limits_{K=1}^{\infty}2^{-K}e^{2\pi i Ky})=\frac{2\sin(2\pi y)}{5-4\cos(2\pi y)}.
\end{eqnarray*}
The direct calculation then shows that this gives precisely the identity (\ref{relation}). We get a linear relation for each higher odd-order derivative.
We will now show that these linear relations characterize the Minkowski question mark function uniquely. This is a consequence of the following
\begin{prop}
Let $\overline{F}(x)$ be a continuous monotone distribution
function: $\overline{F}(x)=0$, $\overline{F}(\infty)=1$. Moreover,
assume $\overline{F}(x)=O(x)$ as $x\rightarrow 0_{+}$,
$1-\overline{F}(x)=O(x^{-1})$ as $x\rightarrow+\infty$. Suppose, for
each odd $r\in\mathbb{N}$ we have $\int_{0}^{\infty}\log^{r}x\d
\overline{F}(x)=0$. Then $\overline{F}(x)+\overline{F}(1/x)=1$.
\end{prop}
 If we put $x=e^{t}$, $F(e^{t})=G(t)$, then this statement is a
standard result in probability theory, claiming that a continuous distribution on $(-\infty,\infty)$ having an exponential decay at both ends whose odd-order moments vanish satisfies $G(t)+G(-t)=1$.
 \subsection{Special value at $s=1$}
 \label{2.4}
Now, multiply (\ref{trig}) by the sum $\sum_{n=1}^{N}\sin(2\pi n
x)$, integrate over $x\in[0,1]$, and take the limit
$N\rightarrow\infty$ using the Lebesgue-Riemann lemma. This implies
\begin{eqnarray*}
\mathfrak{M}(1)=\sum\limits_{n=1}^{\infty}\frac{d_{n}}{n}=\pi\int\limits_{0}^{1}\frac{(?(x)-x)\cdot\sin (2\pi x)}{1-\cos (2 \pi x)}\d x
=\pi\int\limits_{0}^{1}(?(x)-x)\cot(\pi x)\d x
=-0.45595_{+}.
\end{eqnarray*}
To calculate this integral, we used Riemann sums. If the interval is divided into $55.000$ equal subintervals, this gives $5$ correct digits. This, probably, can be pushed to $6$ digits, but not more. However, there exists a far superior method which provides $30$ correct digits; we can push this, with some additional effort, much further. Here is its description.\\

We have the Taylor series for $\pi\cot(\pi x)$:
\begin{eqnarray}
\pi\cot(\pi x)=\frac{1}{x}-
2\sum\limits_{n=1}^{\infty}\zeta(2n)x^{2n-1},\quad |x|<1.\label{cot}
\end{eqnarray}
Further, by a direct calculation,
\begin{eqnarray*}
\int\limits_{0}^{1}(?(x)-x)x^{2n-1}\d x=\frac{1}{2n(2n+1)}-
\frac{m_{2n}}{2n}.
\end{eqnarray*}
Thus,
\begin{eqnarray}
\mathfrak{M}(1)=\int\limits_{0}^{1}(?(x)-x)\frac{\d x}{x}
-\sum\limits_{n=1}^{\infty}\frac{\zeta(2n)}{n(2n+1)}
+\sum\limits_{n=1}^{\infty}\frac{\zeta(2n)\cdot m_{2n}}{n}.\label{ttar}
\end{eqnarray}
Integration term-by-term can be easily justified. In fact, let us write $?(x)-x=x(1-x)r(x)$; then $r(1)=1$, since $?(1-x)=1-?(x)$ behaves like $2^{-1/x}$ at $x=0$, and the Taylor series for $\pi\cot(\pi x)x(1-x)$ converges uniformely in the interval $[0,x_{0}]$ for every $0\leq x_{0}<2$. As is clear from (\ref{cot}), we have
\begin{eqnarray*}
\sum\limits_{n=1}^{\infty}\frac{\zeta(2n)}{n(2n+1)}=
\int\limits_{0}^{1}\Big{(}\frac{1}{x}-\pi\cot(\pi x)\Big{)}(1-x)
\d x=\log2\pi-1;
\end{eqnarray*}
the last value is taken from the tables. Further, it was proved in \cite{alkauskas_i} that
\begin{eqnarray}
\int\limits_{0}^{1}\log x\d?(x)=-2\int\limits_{0}^{1}\log(1+x)\d?(x).\label{log2}
\end{eqnarray}
We note that the proof of this identity essentially employs the functional equations (\ref{distr}). Using this, we obtain
\begin{eqnarray*}
\int\limits_{0}^{1}(?(x)-x)\frac{\d x}{x}=\int\limits_{0}^{1}\log x\d x-\int\limits_{0}^{1}\log x\d ?(x)=
-1+2\int\limits_{0}^{1}\log(1+x)\d ?(x).
\end{eqnarray*}
Further,
\begin{eqnarray*}
\sum\limits_{n=1}^{\infty}\frac{\zeta(2n)\cdot m_{2n}}{n}=
\sum\limits_{n=1}^{\infty}\frac{(\zeta(2n)-1)\cdot m_{2n}}{n}
+\int\limits_{0}^{1}\sum\limits_{n=1}^{\infty}\frac{x^{2n}}{n}\,\d ?(x).
\end{eqnarray*}
Using the symmetry property $?(x)+?(1-x)=1$ and (\ref{log2}), we can simplify the last integral:
\begin{eqnarray*}
\int\limits_{0}^{1}\sum\limits_{n=1}^{\infty}\frac{x^{2n}}{n}\,\d ?(x)&=&\int\limits_{0}^{1}-\log(1-x^2)\d ?(x)
=-\int\limits_{0}^{1}\log(1-x)\d ?(x)-\int\limits_{0}^{1}\log(1+x)\d ?(x)\\
&=&-\int\limits_{0}^{1}\log x\d ?(x)-\int\limits_{0}^{1}\log(1+x)\d ?(x)
=\int\limits_{0}^{1}\log(1+x)\d ?(x).
\end{eqnarray*}
Also, it was proved in \cite{alkauskas_mc} that
\begin{eqnarray*}
\int\limits_{0}^{1}\log(1+x)\d ?(x)=
-\sum\limits_{n=1}^{\infty}\frac{m_{n}}{n\cdot 2^{n}}+\log 2.
\end{eqnarray*}
 One of the results in \cite{alkauskas_gmj} is that the generating function of moments $m_{n}$ satisfies the three-term functional equation, and this, as an application, gives a very good method to calculate these moments with high accuracy. Collecting everything together in (\ref{ttar}), we obtain
\begin{prop}We have the fast converging series for the constant $\mathfrak{M}(1)$, and it is given by
\begin{eqnarray*}
\mathfrak{M}(1)
=\sum\limits_{n=1}^{\infty}\frac{(\zeta(2n)-1)\cdot m_{2n}}{n}
-3\sum\limits_{n=1}^{\infty}\frac{m_{n}}{n\cdot 2^{n}}
-\log\pi+\log 4.
\end{eqnarray*}
\end{prop}
(Compare this with (\ref{evenmo})). As was proved in \cite{alkauskas_lmj}, the constants $m_{n}$ asymptotically behave, up to the constant factor, as $n^{1/4}{\sf C}^{-\sqrt{n}}$, ${\sf C}=e^{-2\sqrt{\log 2}}=0.18916_{+}$. Further, $(\zeta(2n)-1)\sim 2^{-2n}$. Thus, both sums have a good rate of convergence. This gives
\begin{eqnarray*}
\mathfrak{M}(1)=\sum\limits_{n=1}^{\infty}\frac{d_{n}}{n}=-0.455959203740245619075047841829_{+};
\end{eqnarray*}
all digits are exact. Other special values of $\mathfrak{M}(s)$ at odd positive integers can be numerically calculated using analogous techniques.

\end{document}